\documentclass[a4paper,12pt]{amsart}
\usepackage{amsmath,amsthm}
\usepackage[longnamesfirst]{natbib}
\numberwithin{equation}{section}

\newtheorem{thm}{Theorem}[section]
\newtheorem{prop}[thm]{Proposition}
\newtheorem{lem}[thm]{Lemma}

\theoremstyle{definition}

\newtheorem{as}[thm]{Assumption}
\theoremstyle{remark}
\newtheorem{rem}[thm]{Remark}
\def\R{ {\mathbf R} }
\def\C{ {\mathbf C} }

\def\Z{ {\mathbf Z} }

\def\Ex{ \mathbf{E} }

\newcommand{\bE}{\mathbf{E}}
\newcommand{\bN}{\mathbf{N}}
\newcommand{\bP}{\mathbf{P}}
\newcommand{\bR}{\mathbf{R}}

\newcommand{\cD}{\mathcal{D}}
\newcommand{\cX}{\mathcal{X}}
\renewcommand{\P}{\mathbf{P}}

\begin{document}
%\title[Discrete It\^o Calculus applied to He's Framework]
\title{A Discrete It\^o Calculus Approach \\
to He's Framework for \\
Multi-Factor Discrete Markets}
\author{Jir\^o Akahori}
\address{Department of Mathematical Sciences \& Research Center for Finance,
Ritsumeikan University 1-1-1, Nojihigashi, Kusatsu, Shiga, 525-8577, Japan}
\email{akahori@se.ritsumei.ac.jp}
\date{}
\thanks{This research was supported by Open Research Center Project for Private Universities: matching fund subsidy from MEXT, 
2004-2008
and also by Grants-in-Aids for Scientific Research (No. 18540146)
from the Japan Society for Promotion of Sciences.}

\pagestyle{plain}
\maketitle
\begin{abstract}
In the present paper, a discrete version of It\^o's formula for a class of multi-dimensional
random walk is introduced 
and applied to the study of a discrete-time complete market model 
which we call He's framework. %in mathematical finance.
The formula unifies continuous-time and discrete-time settings and
by regarding the latter as the finite difference scheme of the former,
the order of convergence is obtained. 
The result shows 
that He's framework cannot be of order $1$ scheme 
except for the 
one dimensional case.
%which is known as CRR model.
\end{abstract}

\section{Introduction}
In \citet{he:90}, the binomial tree approach by \citet{CRR}
is generalized
to a multi-nomial one and limit theorems concerning
pricing kernels and hedging strategies are established.
The main feature of He's multi-nomial tree framework 
is that each approximating market
itself is arbitrage-free and \underline{\bf complete}. 
%By this reason, 
%one can say that it is continuous-time models that approximate He's discrete ones.

In the present paper, a new insight to He's framework,
which leads to further applications,
will be introduced. 
The insight comes from
a discrete version of It\^o's formula.
%which will be called {\bf Szabados-Fujita formula}.
As is the case with continuous-time models,
our discrete It\^o formula relates the value process
of a contingent claim to a difference equation.
This means that
the formula enables a discrete version 
of so-called {\em partial differential equation (PDE) approach}
to the pricing-hedging problems in the literature of mathematical finance;
we do not use the usual martingale argument.

Further, if a continuous-time limit exists, then the discrete equations obtained 
via our It\^o formula can be seen as {\em explicit} finite difference approximations of
the limit PDE, and we can obtain the order of convergence by using the standard argument of the
finite difference scheme. 

The contributions of the present paper are:
\begin{itemize}
	\item a multi-dimensional version of discrete It\^o formula [Theorem \ref{ItoF}]
which enables the discrete PDE approach.
	\item the order of convergence of the value functions of European %/American 
options within He's framework [Theorem \ref{limittheorem}],
which is proved to be %at least 
$ \mathrm{O} (N^{-1/2}) $ in general and $ \mathrm{O}(N^{-1}) $
in single risky asset cases. Here $ N $ is the number of time-discretization steps.%, Theorem \ref{American1}].
%\item to show that a higher order scheme,
%where the order of convergence of is at least $ \mathrm{O}(N^{-1}) $,
%is within our scope [Proposition \ref{EuroDis}].
%\item to relate He's framework to 
%the representation theory of finite groups [Section \ref{rep}].
\end{itemize}
The point is that 
{\em completeness makes it slow};
as He's framework is based on completeness of market.
The first observation from discrete-It\^o formula 
shows that approximations by the discrete market model 
of He's framework are always a kind of
finite-difference approximation of a PDE, 
while the second observation says that the convergence
is much slower when $ n \geq 2 $ than
an approximation by an Euler-Maruyama scheme of
first order.

\ 

This paper is organized as follows.
In section \ref{Heov}, a quick review of He's framework
will be presented.
In section \ref{Ito}, the Szabados-Fujita formula and the discrete PDE framework will be introduced. 
In section \ref{limit}, a limit theorem will be established. 
%In section \ref{higher}, 
%we will relate a higher order finite difference scheme
%to our discrete It\^o calculus.
In section \ref{rep},
the relations with the group theory will be explained. 
Finally in section \ref{proof},
proofs of the theorems 
in the present paper will be undertaken.

\begin{rem}
This paper is motivated by
the textbook \citet{Fujita_Japanese_text},
where he gives a very nice description of 
{\em from CRR to Black-Scholes} argument
by using his discrete It\^o formula.
His (and our) approach would be very instructive
for those who are not familiar with higher mathematics.
\end{rem}

\noindent {\bf Acknowledgment.} 
The author wishes to acknowledge 
the hospitality extended to him by Professor Marc Yor
during a sabbatical stay at Paris VI, when 
the first version of the present paper was written.
The author also acknowledge suggestions 
by Professors Freddy Delbaen, Raouf Ghomrasni, and
anonymous referees.
%for his/her very useful comments.  

%%%%%%%%
\section{He's framework: an overview}\label{Heov}
In essence, \citet{he:90} approximated 
$ n $-dimensional Brownian motions by a system of
mutually orthogonal martingales of finite states--- 
$ (n + 1) $ states at each step.

Let us briefly review He's framework.
%The problem lies in 
%the choice of an orthonormal basis of 
%$ \R^{n+1} $ which includes the vector of 
Let $ %\mathbf{e} = 
( e_{i,j} )_{0 \leq i,j \leq n} $ 
be an $ (n+1) \times (n+1) $-orthogonal matrix 
such that $ e_{0,j} > 0 $ for $ j=0,1,...,n $, 
and define 
\begin{equation}\label{Eset}
\mathcal{E} 
:= \left\{ {e}_j =  \frac{1}{e_{0,j} }(e_{1,j},...,e_{n,j} ) \in \R^n \,:\,
j=0,...,n \right\}.
\end{equation}
Let $ \tau \equiv (\tau^1,...,\tau^n) $ 
be a random variable taking values
in $ \mathcal{E} $ with
\begin{equation*}
\P (\tau = {e}_j ) = e_{0,j}^2, \,\,j=0,1,...,n.
\end{equation*}
Then, \footnote{Note that the converse is true;
any random variable $ \tau $ satisfying (\ref{cov})
is, if it is defined on a finite set,
constructed in the above way 
from such an orthogonal matrix.
\citet{he:90} treated only the uniform cases
of $ e_{0,0}= \cdots = e_{0,n} = 1/\sqrt{n+1} $.}
\begin{equation}\label{cov}
\Ex [ \tau^i ] = 0, \,i=1,...,n, \quad
\mbox{and} \quad  \mathrm{Cov} (\tau^i, \tau^j) 
= \begin{cases}
1 & (i=j) \\
0 & (i \ne j ).
\end{cases}
\end{equation}

Let $ \tau_1,....,\tau_t,... $ be independent copies of $ \tau $.
Define a sequence of $ \R^n $ valued stochastic processes $ \{ X^{N}\} $ by
\begin{equation*}
X^{N}_t = X_0 + N^{-1/2} \sum_{u=1}^{[Nt]} \tau_u
\end{equation*}
for a given initial point $ X_0 \in \R^n $. 
By (\ref{cov}), 
components of $ X^N_t - X^N_0 $
are mutually orthogonal martingales, 
and therefore, the martingale central limit theorem 
(see \citet{MR88a:60130} for example) ensures that the law of $ X^N $ 
converges weakly to the $ n $-dimensional Wiener measure as $ N \to \infty $.

Fix $ N \in \bN $. 
For $ T>0 $, we denote $ T_N=[TN]/N $. 
For a subinterval $ I $ of $ [0,\infty ) $, we denote 
$ I_N = I \cap \{ k/N : k=0,1,2,\ldots \} $.

In our market there are $ (n+1) $-securities 
whose prices are given by
$ S^{j,N}_t \equiv  h^{j,N} ( t, X^N_t )$
for $ j=0,1,...,n $, where 
$ h^{j,N} $'s are real functions defined on 
$ [0,T]_N \times \R^n $ such that 
the following $ (n+1)\times(n+1) $-matrix 
\begin{equation*}
H^N (t,x) := 
\begin{pmatrix}
h^{0,N} (t, x+ N^{-1/2}e_0 ) & \cdots & h^{n,N} (t, x+N^{-1/2}e_{0} ) \\
h^{0,N} (t, x+ N^{-1/2}e_1 ) & \cdots & h^{n,N} (t, x+N^{-1/2}e_{1} )  \\
\vdots & \ddots & \vdots \\
h^{0,N} (t, x+ N^{-1/2}e_n ) & \cdots & h^{n,N} (t, x+N^{-1/2}e_{n} ) 
\end{pmatrix}
\end{equation*}
is invertible for arbitrary $ (t,x) \in [0,T]_N \times \R^n $. 

Suppose that at time $ t=k/N $ we have 
$ \theta^j_t \equiv \theta^j_t (\tau_1,...,\tau_k) $ 
amount of $ j $-th security for each $ j \in \{0,1,...,n \} $.
The cost of the portfolio at time $ t $ is 
\begin{equation}\label{cost}
v^N (t, \tau_1,...,\tau_k) := \sum_j h^{j,N} (t,X^N_t) \theta^j_t, 
\end{equation}
and at time $ t + N^{-1} $ the value of the 
portfolio becomes 
\begin{equation*}
v^N (t+N^{-1}, \tau_1,...,\tau_k, \tau_{k+1}) :=
\sum_j h^{j,N} (t,X^N_{t + N^{-1}}) \theta^j_t ,
\end{equation*}
or equivalently
\begin{equation}\label{value}
\begin{pmatrix}
v^{N} (t+N^{-1},\tau_1,...,\tau_k, N^{-1/2}e_0 ) \\
v^{N} (t+N^{-1},\tau_1,...,\tau_k, N^{-1/2}e_1 ) \\
\vdots  \\
v^{N} (t+N^{-1},\tau_1,...,\tau_k, N^{-1/2}e_n ) 
\end{pmatrix}
= H^N (t+ N^{-1},x) 
\begin{pmatrix}
\theta^0_t \\
\theta^1_t \\
\vdots \\
\theta^n_t
\end{pmatrix}.
\end{equation}
If the portfolio is self-financed, 
then $ c^N (t,\cdot) = v^N (t, \cdot) $.
Since we have assumed that $ H^N $ is invertible, 
we have by combining (\ref{cost}) and (\ref{value}),
\begin{equation}\label{recursivep}
\begin{split}
%&v^N (T_N,x) = \Phi^N (x); \quad x \in \R^n, \\
&v^N (t,x) = \\
&(h^{0,N} (t, x),...,h^{n,N} (t, x) )
H^N (t+ N^{-1},x)^{-1} %\\&\times 
\begin{pmatrix}
v^{N} (t+N^{-1},x,N^{-1/2}e_0 ) \\
v^{N} (t+N^{-1},x,N^{-1/2}e_1 ) \\
\vdots  \\
v^{N} (t+N^{-1},x,N^{-1/2}e_n ) 
\end{pmatrix}
; \\
& \quad t \in [0,T_N)_N, \quad x  \in \mathcal{E}^{[tN]}.
\end{split}
\end{equation}
If the terminal value (to be hedged) 
$ \Phi : \mathcal{E}^N \to \R $ is 
dependent only on $ X^N_T $, 
then $  v^N (T-N^{-1}, \cdot ) $ depends only on 
$ X^N_{t - N^{-1}} $ etc, etc, and finally 
we have the following recursive equation,
which has a unique solution:  %(\ref{recursiveE}).
\begin{equation}\label{recursiveE}
\begin{split}
&v^N (T_N,x) = \Phi^N (x); \quad x \in \R^n, \\
&v^N (t,x) = \\
&(h^{0,N} (t, x),...,h^{n,N} (t, x) )
H^N (t+ N^{-1},x)^{-1} %\\&\times 
\begin{pmatrix}
v^{N} (t+N^{-1},x+N^{-1/2}e_0 ) \\
v^{N} (t+N^{-1},x+N^{-1/2}e_1 ) \\
\vdots  \\
v^{N} (t+N^{-1},x+N^{-1/2}e_n ) 
\end{pmatrix}
; \\
& \quad t \in [0,T_N)_N, \quad x  \in \R^n.
\end{split}
\end{equation}
Here all we can say is that  
$ v (t,X^N_t) $ is the replication cost 
(at time $ t $) of 
an European option whose pay-off is described by 
$ \Phi (X^N_T) $, where $ \Phi^N : \R^{n+1} \to \R $.

As is well known, 
absence of arbitrage opportunities is equivalent 
to the positivity of the state price (see e.g. \citet{duffie:96}).
In other words, 
denoting $ h^N (t,x) = $ 
$(h^{0,N} (t, x) $
$,...,$
$ h^{n,N} (t, x) ) $, 
\begin{equation}\label{NA}
\text{each component of $ h^N (t,x)H^N (t+N^{-1},x)^{-1} $ is strictly positive.}
\end{equation}
Under the hypothesis of (\ref{NA}), the unique solution 
$ v^N(t,x) $ is the unique fair price 
at time $ t \in [0, T) $ and the state $ x \in \R^n $
with $ S^{j,N}_t =  h^{j,N} ( t, x ) $
of the European option whose pay-off is 
$ \Phi^N (X^N_T) $.
Note that the invertibility of $ H $ 
is equivalent to completeness of the market.

\ 

The above derivation of (\ref{recursiveE})
is also valid for 
%$ X^N $ is replaced with 
any Markov process\footnote{In general it 
is represented by some $ F_j, j=0,1,...,n $ as 
\begin{equation*}
Z^N_{t + N^{-1}} - Z^N_t = \sum_{j=0}^n
F_j (Z_t^N) \tau_{t+N^{-1}}^j,
\end{equation*}
with a convention of $ \tau^0 \equiv 1 $.
This is because $ 1,\tau^1,...,\tau^n $ forms 
an orthonormal basis of the space of random variables
on generated by $ \tau $. In 
particular, a discrete approximation 
of an SDE by a Markov chain 
always has an Euler-Maruyama representation.}
$ Z^N $ replacing $ X^N $. 
In fact \citet{he:90} modeled the price vector
$ \mathbf{S}_t=(S^1_t,...,S^n_t) $ directly
(meaning $ h^j $'s are identity maps)
by an Euler-Maruyama approximation
of a stochastic differential equation.
Here we have changed the setting as above.
The differences is that we have preserved 
the structure of so-called recombining tree:
if we consider $ \mathbf{S}_{k/N} $ as a function of
$ \tau_1,...,\tau_k $, we have
\begin{equation}\label{recom}
\mathbf{S}_{k/N} (e_{i_1},...,e_{i_k} )
= \mathbf{S}_{k/N} (e_{i_{\sigma(1)}},...,e_{i_{\sigma(k)}} )
\end{equation}
for arbitrary permutation $ \sigma \in \mathfrak{S}_k $.

The reasons for this modification are:
(i) Euler-Maruyama approximations by 
finite-points random variables 
using Monte-Carlo are not practical, (ii)
nor is solving an equation like (\ref{recursiveE}) 
without recombining structure of (\ref{recom}).

In fact, it relaxes quite a lot computational complexity,
by which we mean {\em how many times we need to 
solve the one-step linear equation (\ref{recursiveE})
to obtain the value for $ v^N (t, x) $.} 
In other words, it is the number $ \sharp \cX(t,x, Z) $
of the possible states 
\begin{equation*}
\cX (t,x,Z):= \{ y \in \R^n : \P ( Z^N_{t}  = x, X^N_T =y ) > 0  \}.
\end{equation*}
In general   
we have $ \sharp \cX(T-kN^{-1}, x, Z ) = (n+1)^k $.
Even if $ Z $ is an Euler-Maruyama approximation of
a solution to SDE, almost always this is the case. 
However, the symmetry (\ref{recom}),  
which comes from that of $ X_t $,
reduces it dramatically.
More precisely, we have the following.
\begin{prop}\label{tree}
\begin{equation*}
%\sharp \cX(T-kN^{-1}, x, \mathbf{S}^N )\leq 
\sharp \cX(T-kN^{-1}, x, X^N )= \frac{(k+n)!}{k!n!}.
\end{equation*}
\end{prop}
\begin{proof}
Since $ \{ e_1,...,e_n \} $ spans 
$ n $-dimensional subspace in $\R^{n+1}$,
they have no linear dependence other than $ e_1 + \cdots + e_n = 0 $.
Therefore, the number is equal to
that of solutions to 
\begin{equation*}
x_1 + x_2 + \cdots + x_{n+1} = k, \quad x_j \in \Z_+,\,j=1,...,n+1,
\end{equation*}
which is exactly $ (n+k)!/k!n! $. 
\end{proof}
%%%
\begin{rem}\label{NAcomp}
Denoting by $ A (t,x) $ 
the sum of all the components of
$ h^N (t,x) H^N (t+N^{-1},x)^{-1} $, 
the value process of money market account is given by
\begin{equation*}
\prod_k^{[Nt]} 1/A(k/N,X^N_{(k-1)/N}).
\end{equation*}
In particular,
{\em positive interest rate} 
is equivalent to $ A (t,x) < 1 $ for arbitrary $ (t,x) $.
\end{rem}

%%%
%%%
\section{A Discrete It\^o formula and discrete PDE}\label{Ito}
%To prove the above Theorem \ref{European1}, %and \ref{American1},
Let us introduce a discrete version of It\^o's formula
for the process $ X^N \equiv (X^{N,1},...,X^{N,n}) $. 
\begin{thm}\label{ItoF}%[It\^o-Szabados-Fujita formula]
{\em (i)} For a function $ f : [0,\infty) \times \R^n \to \R $,
we have
\begin{equation}\label{itofujita}
\begin{split}
&f(t,X^N_t) - f(0,X_0) \\
& \qquad = \sum_{u=1}^{[Nt]} \bigg( \sum_{k=1}^n 
\partial^N_k f (u/N, X^N_{(u-1)/N})\,
(X^{N,k}_{u/N} - X^{N,k}_{(u-1)/N}) \\
& \qquad \qquad \quad + ( \frac{1}{2}\Delta^N  + \partial_t^N )
f(u/N, X^{N}_{(u-1)/N}) /N \bigg),
\end{split}
\end{equation}
where
\begin{equation}\label{diff}
\begin{split}
\partial_{k}^N f (\cdot, x) &= \sqrt{N}
\sum_{j=0}^{n} f( \cdot , x + N^{-1/2} e_{j} )e_{0,j} e_{k,j}  \\
\Delta^N f (\cdot,x) &= 2 {N} 
\sum_{j=0}^{n} 
\{ f( \cdot , x + N^{-1/2} e_j) - f (\cdot, x) \} e_{0,j}^2, \\
\partial^N_t f( t, \cdot ) &= N(f( t, \cdot ) - f(t-N^{-1}, \cdot)).
\end{split}
\end{equation}
{\em (ii)} If $ f $ is in $ C^{1,2} $ in a neighborhood of $ (t,x) $, then
letting $ N \to \infty $, we have
\begin{equation}\label{point}
\partial_j^{N} f (t, x) \to \frac{\partial}{\partial x_j}f(t, x),\,
\Delta^N f (t,x) \to \Delta f (t,x), \partial_t^N f (t,x) 
\to \frac{\partial}{\partial t} f (t, x).
\end{equation}
Here $ \Delta $ is the Laplacian in $ \R^n $.
{\em (iii)} Further, for fixed $ t \in [0,T] $,
if $ f(t,\cdot)$ is in $ C^{3} $ in an open set $ U \subset \R^n $,
then for every compact subset $ K \subset U $, there exists a positive constant $ C_K $
depending only on $ f (t,\cdot) $
such that
\begin{equation}\label{orders}
\max_j \left| \partial_j^{N} f (t, x) - \frac{\partial}{\partial x_j}f(t, x) \right|
+ \left| \Delta^N f (t,x) - \Delta f(t,x) \right|
\leq C_K N^{-1/2}  
\end{equation}
for all $ x $ in $ K $. 
(iv) The order of convergence cannot be improved 
for general $ f \in C^{1,4} $ when $ n \geq 2 $.
(v) For the case of $ n =1 $, 
it %the right-hand-side of (\ref{orders})
can be improved to be $ N^{-1} $, provided that $ f \in C^{1,4} $.
\end{thm}
A proof of Theorem \ref{ItoF} will be given in section \ref{PDIF}.

\begin{rem}
This version of It\^o's formula is different
from those for jump semimartingales which, for example, is appearing in \citet{MR2020294},
different in that ours gives the Doob decomposition of $ f (t,X_t) $.
This version of discrete It\^o's formula 
was introduced by \citet{Fujita_formula} 
for the case of $ n=1 $. 
\citet{MR684465} and \citet{MR92i:60105} also 
studied discrete It\^o formulas as discrete-analogues of 
the standard one, which point of view is what we share in this paper.
It is true that 
it should be called Kudzhma-Szabados-Fujita formula, 
but here the term {\em discrete It\^o formula} 
is preferred since the true name is too long and confusing
\end{rem}

We claim that the recursive equation
(\ref{recursiveE}) defines a discrete PDE 
with respect to these {\em differentials} of (\ref{diff}).
Define 
\begin{equation*}
\Sigma^N (t,x) := \begin{pmatrix}
h^{0,N} (t-N^{-1},x) & \cdots & h^{n,N} (t-N^{-1},x) \\
\partial^N_1 h^{0,N} (t,x) & \cdots & \partial^N_1 h^{n,N} (t,x) \\
\vdots & \ddots & \vdots \\
\partial^N_n h^{0,N} (t,x) & \cdots & \partial^N_n h^{n,N} (t,x)
\end{pmatrix}.
\end{equation*}

\begin{thm}\label{Euro}
Let us assume that the market is arbitrage-free and complete.
Namely, the existence of $ (H^N)^{-1} $ and (\ref{NA}) are assumed.
Then, $ \Sigma^N $ is always invertible 
and $ v^N $ satisfies the following discrete PDE.
\begin{equation}\label{disdif}
\begin{split}
& \nu^N(T_N,x) = \Phi^N(x) ; \quad x \in \bR^n, \\
& \partial ^N_t \nu^N + \frac{1}{2} \Delta^N \nu^N 
- \langle b^N,\nabla^N \nu^N \rangle - c^N (1^N \nu^N)= 0 ; \\
& \qquad t \in (0,T_N)_N , \ x \in \bR^n . 
\end{split}
\end{equation}
Here $ 1^N \nu^N(t,x) = \nu^N(t-N^{-1},x) $ 
and $ (c^N,b^N)=(\partial ^N_t h^N 
+ \frac{1}{2} \Delta^N h^N) [\Sigma^N]^{-1}$. 
\end{thm}
A proof of Theorem \ref{Euro} will be given in section \ref{PEuro}.

\ 

The equation (\ref{disdif}) can be obtained directly by using
the discrete It\^o's formula (\ref{itofujita})
if we a priori assume that $ \Sigma^N $ is invertible.
Let us write $ d Y_t :=  Y_t - Y_{t-N^{-1}} $ for a process $ Y $, 
$ dt := 1/N $, $ \nabla^N = (\partial^N_1,...,\partial^N_n ) $,
$ V_t := v^N (t,X^N_t) $, and so on.
If we have $ d V_t = \sum_{j=1}^{n+1} \theta^j_t d h_t^{j,N} $ and 
$ V_t = \sum_{j=1}^{n+1} \theta^j_{t+1} h^{j,N}_t $, 
then $ \theta $ is the hedging strategy and the problem is settled.
This can be done quite easily in a parallel way with 
the continuous-time cases. 
In fact, we have
\begin{equation*}
\begin{split}
dV &=\nabla v^N \cdot dX^N + \left( \partial^N_t v^N + \frac{1}{2} 
\Delta^N v^N \right) \,dt = \sum_{j=1}^{n+1} \theta^j d h^{j,N}, \\
dh^{j,N} &= \nabla h^{j,N} \cdot dX^N 
+ \left( \partial^N_t h^{j,N} + \frac{1}{2} \Delta^N h^{j,N} \right) \,dt,
\end{split}
\end{equation*}
and $ v^N = \sum_{j=1}^{n+1} \theta^j_{t+1} h^{j,N}_t $, hence
$ \theta= (\theta^1,...,\theta^{n+1}) = (\Sigma^N)^{-1} (v^N (t-N^{-1}), \nabla v^N ) $.

Note that the above argument can be applied to the case of $ N = \infty $,
where $ X^\infty $ is the standard Brownian motion, 
$ dt $ is the standard one,  
and so on. 
The corresponding standard PDE shares the algebraic structure with
the discrete ones. 
Since the second assertion (ii) of the above theorem \ref{ItoF}
can be seen as {\em consistency} of the difference operators
of (\ref{diff}) in the context of
finite difference method (see e.g. \citet{MR95b:65003}) , $ v^N $ converges to a solution $ v^\infty $ to the PDE
at least when $ v^\infty $ is regular enough. We will make a detailed study 
about this topic in the next section.

%%%%%%%%%%%%
\section{Limit Theorem}\label{limit}
The solution $ {v^N(t,\cdot):\bR^n \to \bR} $ 
is solved inductively for each $ t \in {[0,T_N)_N} $, 
and for each $ {x \in \bR^n} $, 
the function $ {v^N(\cdot,x)} $ on $ {[0,T]_N} $ 
can be extended to a piecewise-constant function on $ {[0,T]} $. 
We choose such an extension on $ [0,T] \times \R^n $
and denote it by the same symbol.

Here we assume the followings to establish
our limit theorem.
\begin{as}\label{standing}
(i) The market is arbitrage-free and complete; i.e. 
we assume (\ref{NA}) and invertibility of $ H^N $. 
(ii) There exist bounded measurable functions
$ b : [0,T] \times \R^n \to \R^{n} $ and $ \Phi : \R^n \to \R $ 
and a continuous function $ c : [0,T] \times \R^n \to \R $
such that
\begin{equation}\label{orderorder}
\sup_N \sup_{t \in [0,T]_N , \ x \in \bR^n} 
N^{1/2}
\left( | b - b^N | + | c - c^N |
+ |\Phi - \Phi^N| \right) < \infty . 
\end{equation}
and 
(iii) they are regular enough to allow
the following partial differential equation to have a bounded solution
in $ C^{1,3} $ whose first order derivatives are also bounded.
\begin{equation}\label{EuroPDE}
\frac{\partial v}{\partial t}
+ \frac{1}{2}\Delta v - \langle b(t,x), \nabla v \rangle - c (t,x) v = 0, 
\quad t \in [0, T), 
\quad v(T, x) = \Phi (x), \quad 
\end{equation}
where $ \Delta $ is the Laplacian of $ \R^n $ and
$ \nabla $ is the gradient operator in $ \R^n $. 
(iv) We also assume that the interest rate is positive 
(See Remark \ref{NAcomp}).
\end{as}
As we pointed out in Remark \ref{NAcomp}, the assumption (iv) 
is equivalent to $ A < 1 $,
and hence, as we shall show in the proof of Theorem \ref{Euro}, 
is equivalent
to the positivity of the first component of 
$ (\partial^N_t h^N + \frac{1}{2} \Delta^N h^N ) [\Sigma^N]^{-1} $.
This in turn implies that $ c $ is positive.

Under the Assumption \ref{standing}, we have the following
\begin{thm}\label{limittheorem}
The solutions $ v^N $ to (\ref{disdif}) converges uniformly
on compact intervals of $ [0,T] \times \R^n $ to the solution 
$ v \in C^{1,3} $ to (\ref{EuroPDE}) in an order of 
$ N^{-1/2} $.
For the one dimensional case, the order can be improved to be $ N^{-1} $
provided that $ v $ is in $ C^{1,4} $ and the order in (\ref{orderorder}) is
replaced with $ N^{-1} $.
\end{thm}
A proof of Theorem \ref{limittheorem} 
will be given in section \ref{PLMT}.

\begin{rem}. 
Our scope covers as a special case the Black-Scholes economy by setting  
$ h^{j,N} (t,x) \equiv S^j_0 e^{ \langle \sigma_j ,x \rangle - \mu t} $ for $ j=1,...,n $
and $ h^{n+1,N}(t,x) \equiv e^{rt} $, where $ S^j_0,  r \in \R_+ $, $ \mu^j \in \R $
and $ \Sigma \equiv [ \sigma_1,...,\sigma_n ] $ is a $ n \times n $ positive definite matrix. 
%Further, as is shown in \citet{AAN}, 
%the framework is almost order-made for 
%discrete-time interest rate models.
\end{rem}

\begin{rem}
Since we are working on a reference measure which is not
necessarily a risk neutral measure nor so-called {\em physical measure},
$ W $ can be a diffusion process other than Brownian motions
under those measures.
Roughly speaking, $ W $ can be a solution (in the {\em weak} sense!) to a stochastic 
differential equation whose diffusion coefficients are constant functions.
In one dimensional cases, by scaling we can work on any diffusion whose
diffusion coefficient is monotone and smooth.
\end{rem}

%%%%%%%%%%%%

%%%
\section{Supplementary Remark: 
relations to Group Representation}\label{rep}

%As mentioned in Remark \ref{remgroup}, 
We remark here that a specification of $ \mathcal{E} $ 
can be done 
with the help of group representation theory.

Let us recall the basics of group representation theory (see e.g. \citet{MR80f:20001}).
Let $ G $ be a compact abelian group, and let $ \widehat{G} $ be
its dual group. The members of $ \widehat{G} $ are often called {\em characters},
which forms an orthonormal basis of $ L^2(G; \C) $; the space of square integrable
functions on $ G $ over $ \C $ with respect to its Haar measure.
Since $ L^2 (G;\C) $ is a complex vector space,
we need to modify it to get an orthonormal basis over real field.
One candidate is obtained by the transform $ \varphi : \C \to \R $ defined by $ \varphi (x + i y) = x + y $.
It is easy to check that $ \{\varphi(\chi): \chi \in \widehat{G} \} $
is a orthogonal basis of $ L^2(G; \R) $.
The group $ \widehat{G} $ always contains a unit,
which corresponds to $ \mathbf{1} $.
Thanks to Peter-Wyel Theorem, the above argument is extended to
non-abelian groups.

In particular, a choice of group $ G $ with $ |G| = n+1 $
gives us an $ \mathcal{E} $. 
The simplest choice may be the cyclic group $ C_{n+1} $.
In this case, $ \tau = (\tau^1,....,\tau^{n}) $ is 
obtained by taking $ \tau^k = \varphi(\eta^k) $
where $ \eta $ is a uniformly distributed 
random variable taking values in $ (n+1) $-th units of root. 
The fundamental theorem of finitely generated abelian groups
says that the characters are always taking values in a set of the units of root. 
Therefore, the {\em scenarios} are generated by a random walk on 
a ring of {\em integers} of an algebraic number field.
The easiest case ($n=2$) is studied in \citet{Aka_Parisian}. 
%For details, see \citet{AAN}. %and \citet{Aka_Parisian}.
%\end{rem}

%%%
\section{Proofs}\label{proof}

\subsection{A Proof of Theorem \ref{ItoF}}\label{PDIF}
Let $ L (\tau) $ be a linear space of $\tau $-measurable 
real valued random variables.
Since $\tau $ takes only $ (n+1) $-distinct values, 
the dimension of $ L (\tau) $ is $ (n+1) $. 
On the other hand,
as a matter of course, the coordinate maps 
$ \tau^1,...,\tau^n $ are members of $ L (\tau) $.
The moment condition (\ref{cov}) 
says that $ \{ \tau^1,....,\tau^n \} $ and 
constant function $ \mathbf{1} $ are mutually orthogonal 
with respect to the inner product 
$ \langle x,y \rangle =\Ex [xy] $.
Hence $ \{ \mathbf{1}, \tau^1,...,\tau^n \} $
is an orthonormal basis of $ L (\tau) $.

Orthogonal expansion 
of $ f (t, x + N^{-1/2} \tau) $ with respect to the basis 
 $ \{ \mathbf{1}, \tau^1,...,\tau^n \} $ are as follows:
\begin{equation*}\label{Fourier2}
\begin{split}
f (t, x + N^{-1/2} \tau)
&= \sum_{k=1}^n \Ex [f (t, x + N^{-1/2} \tau)\tau^k ] \tau^k
+ \Ex [f (t, x + N^{-1/2} \tau)] \\
&= \sum_{k=1}^n 
\left( \sum_{j=0}^{n} \P (\tau=e_j) f(t,x+ N^{-1/2} e_j ) 
\frac{e_{k,j}}{e_{0,j}} \right) \tau^k  \\
& \hspace{3cm}
+ \sum_{j=0}^{n} \P (\tau=e_j) f(t,x+ N^{-1/2} e_j ) \\
&= \sum_{k=1}^n 
\left( \sum_{j=0}^{n} f(t,x+ N^{-1/2} e_j ) 
e_{0,j} e_{k,j} \right) \tau^k  \\
& \hspace{3cm}
+ \sum_{j=0}^{n} e_{0,j}^2 f(t,x+ N^{-1/2} e_j ) \\
&= \sum_{k=1}^n \partial^N_k f (t,x) \frac{\tau^k}{\sqrt{N}}
+ \frac{1}{2N} \Delta^N f (t,x)\, 
+ f(t,x).
\end{split}
\end{equation*}

Substituting $ X^{N}_u $ for  $ x $
and $ X^{N,k}_{u+N^{-1}} - X^{N,k}_u $ for $ \tau^k/\sqrt{N} $, 
we have
\begin{equation}\label{ItoFD}
\begin{split}
&f(u+ N^{-1}, X^N_{u + N^{-1}})-f(u,X_u^N) \\
&\,\, 
= \sum_{k=1}^n \partial^N_k f (u+ N^{-1}, X^N_u) ( X^{N,k}_{u+N^{-1}} - X^{N,k}_u ) \\
& \qquad
+ ( \frac{1}{2} \Delta^N + \partial^N_t ) f (u+ N^{-1}, X^N_u)/N.
\end{split}
\end{equation}
By summing up (\ref{ItoFD}) for $ u=0, 1/N, 2/N,...,([Nt]-1)/N $, 
we obtain (\ref{itofujita}).

%(ii), (iii)
Let us consider next the following formal Taylor expansion of 
$ f (t, x + N^{-1/2} \tau ) $ with respect to $ N^{-1/2} \tau $:
\begin{equation*}
\begin{split}
&f ( t, x + N^{-1/2} \tau ) \\
&\quad = f(t,x) + 
\frac{1}{\sqrt{N}} \langle \nabla f, \tau \rangle 
+ \frac{1}{\sqrt{N}} \langle \nabla f \otimes \nabla f, 
\tau \otimes \tau \rangle \\
& \hspace{3cm} + \cdots + \frac{1}{N^{m/2}}
\langle (\nabla f)^{\otimes m \,\scriptsize{\mbox{times}} }, 
\tau^{\otimes m \,\scriptsize{\mbox{times}}} \rangle
+ \cdots 
\end{split}
\end{equation*}
Recalling (or observing the proof given above) 
that 
\begin{equation*}
\partial_k^N f (t,x) 
= \sqrt{N} \Ex [f (t, x + N^{-1/2} \tau)\tau^k ],
\end{equation*}
and
\begin{equation*}
\Delta^N f (t,x) = 2 N \Ex [f (t, x + N^{-1/2} \tau)
-f (t, x)],
\end{equation*}
we have the following formal expansions:
\begin{equation*}
\begin{split}
\partial_k^N f (t,x) & =  
\sqrt{N}\,\Ex [f (t, x )\tau^k ] 
+ \Ex [\langle \nabla f, \tau \rangle \tau^k ] \\
& \quad + \frac{1}{2\sqrt{N}} \,\Ex
[\langle \nabla f \otimes \nabla f, 
\tau \otimes \tau \rangle \tau^k ]  \\
& \quad + \frac{1}{2 N} \Ex [\langle \nabla f 
\otimes \nabla f \otimes \nabla f, 
\tau \otimes \tau \otimes \tau \rangle \tau^k  ] + \cdots + \cdots
\end{split}
\end{equation*}
and
\begin{equation*}
\begin{split}
\frac{1}{2}\Delta^N f (t,x) & =  
+ \sqrt{N}\Ex [\langle \nabla f, \tau \rangle ]
+ \frac{1}{2} \,\Ex
[\langle \nabla f \otimes \nabla f, 
\tau \otimes \tau \rangle ]  \\
& \quad + \frac{1}{6\sqrt{N}} \Ex [\langle \nabla f 
\otimes \nabla f \otimes \nabla f, 
\tau \otimes \tau \otimes \tau \rangle] \\
& \quad + \frac{1}{24 N} \Ex [\langle (\nabla f)^{\otimes 4}, 
\tau^{ \otimes 4} \rangle ]
+ \cdots + \cdots.
\end{split}
\end{equation*}
Now the assertions (ii) and (iii) are verified 
since for $ f(t,\cdot) \in C^k $ the expansion up to $ k $-th 
term is valid.

The assertion (v) is verified by looking at the case
$ \P (\tau = \pm 1 )= 1/2 $ where $ \Ex [\tau^3 ] = 0 $. 
The assertions (iv) is a consequence of the following lemma
\qed 

\begin{lem}
Suppose that $ \tau: \Omega \to \R^n $ satisfies 
(\ref{cov}) and $ \sharp \Omega = n+1 $. Then 
if $ n \geq 2 $, there exists $ (i,j,k) \in \{ 1,...,n \}^3 $
such that $ \Ex [\tau^i \tau^j \tau^k ] \ne 0 $.
\end{lem}
\begin{proof}
Denote $ \mathcal{D} = (e_{i,j})_{0 \leq i,j \leq n} $
and let 
\begin{equation*}
\mathcal{D}_k = \mathrm{diag} 
[e_{k,0}/e_{0,0}, e_{k,1}/e_{0,1}, \cdots, e_{k,n}/e_{0,n}],
\,\,k=1,...,n.
\end{equation*}
Then 
one will find that the $ (i,j) $-th component of 
$  \cD^* \cD_k \cD  $ is given by 
\begin{equation*}
d_{i,j} = \sum_{l=0}^n \frac{e_{i,l}e_{k,l} e_{j,l}}{e_{0,l}}
= \Ex [\tau^i \tau^k \tau^j],
\end{equation*}
where conventionally $ \tau^0 \equiv 1 $ and
$ d_{i,j} $ for $ 0 \leq i,j \leq n $ are numbered as follows:
\begin{equation*}
\cD^* \cD_k \cD =
\begin{pmatrix}
d_{0,0} & d_{0,1} & \cdots & d_{0,n} \\
d_{1,0} & d_{1,1} & \cdots & d_{1,n} \\
\vdots & \ddots & & \vdots \\
d_{n,0} & d_{n,1}  & \cdots & d_{n,n} 
\end{pmatrix}
.
\end{equation*}

If we assume
$ \Ex [\tau^i \tau^j \tau^k ] = 0 $
for all $ (i,j,k) \in \{ 1,...,n \}^3 $, 
then 
for arbitrary fixed $ k $
we have $ d_{i,j} = 0 $ for every $ (i,j) \in \{ 1,...,n \}^2 $.
Since $ d_{0,j} = \Ex [ \tau^k \tau^j ] = \delta_{k,j} $,
we notice that 
$ \mathrm{rank} \cD^* \cD_k \cD = \mathrm{rank} \cD_k 
= 2 $. 
This implies, since $ \cD_k $ is a diagonal matrix,
$ e_{k,j} = 0 $ except for exactly two $j$'s, 
for which we write $ k_+ $ and $ k_- $.

We may assume without loss of generality
$ e_{k,k_-} < 0 < e_{k,k_+} $
since $ \Ex[\tau^k] 
=e_{k,k_-} e_{0,k_-} + e_{k,k_+} e_{0,k_+} $ must be 
zero. This in turn implies 
$ \{k_-,k_+\} $, $ k=1,...,n $ must be disjoint
to fulfill $ \Ex[ \tau^k \tau^{k'} ] = 0 $ for $ k \ne k' $. 
Hence finally we notice that
$ 2 n \leq n+1 $. This implies $ n =1 $. 
\end{proof}

%%%
\subsection{Proof of Theorem \ref{Euro}}\label{PEuro}
We will write 
\begin{equation*}
\tilde{f} (x) = ( f(x+ N^{-1/2} {e}_0),...,
f(x+N^{-1/2} e_{n}) ) 
\end{equation*}
for $ f : \R^n \to \R $. 
Note that $ \tilde{f} $ is a map to $ \R^{n+1} $.

As in the above proof we denote
$ \mathcal{D} = (e_{i,j})_{0 \leq i,j \leq n} $.
Then we have 
\begin{equation*}
\mathcal{D} \tilde{f}(x) = \left( f(x) + (2N)^{-1} \Delta^N f(x), 
N^{-1/2} \partial^N_1 f(x),...,
N^{-1/2} \partial^N_n f(x) \right). 
\end{equation*}
Since $ \cD H^N = \cD \tilde{h}^N = \Sigma^N + (a,0,\ldots,0)^* 
$ for some $ a=a(t,x) $, we have
%Since $ [H^N (t,x)]^{-1} \tilde{v^N} (t,x)
%= [ \mathcal{D} H^N (t,x)]^{-1} E \tilde{v^N} (t,x) $, we have
\begin{equation}\label{DM}
\Sigma^N [\mathcal{D} H^N (t,x)]^{-1}
= 
\begin{pmatrix}
\pi^N_1 (t,x) & \cdots &  \pi^N_{n+1} (t,x) \\
\mathbf{0} & \quad N^{1/2} I_n & \,
\end{pmatrix}
\end{equation}
where $ \mathbf{0} = (0,...,0)^* \in \R^n $,
$ I_n $ is the unit $ n \times n $ matrix, and 
\begin{equation}
\pi^N = (\pi^N_1,\ldots,\pi^N_{n+1}) 
= (1^N h^N) [\cD H^N]^{-1} = (1^N h^N)[H^N]^{-1} \cD^{-1} . 
\end{equation}
Since $ \cD^{-1}=\cD^* $, we have 
$ \pi^N_1=1^N A $, 
the sum of the components of 
$ (1^N h^N) [H^N]^{-1} $, 
which is strictly positive by the assumption (i), 
and hence $ \Sigma^N $ is invertible. %} \\

Using (\ref{DM}), we have
\begin{equation}
\begin{split}
\Sigma^N [\cD H^N]^{-1} \cD \tilde{\nu}^N
=& [\pi^N \cD \tilde{\nu}^N,\partial^N_1,\ldots,\partial^N_n]^* \\
=& (1^N,\nabla^N)^* \nu^N. 
\end{split}
\end{equation}
Here we use the equality
$ (1^N h^N) [H^N]^{-1} \tilde{\nu}^N = 1^N \nu^N $ 
which comes from (\ref{recursiveE}). 
Hence we have
\begin{equation}\label{7.4}
\cD \tilde{\nu}^N 
= \cD H^N [\Sigma^N]^{-1} (1^N,\nabla^N)^* \nu^N. 
%\tag{7.4}
\end{equation}
In particular, we have the following relation 
from the first component of the above (\ref{7.4}): 
\begin{equation}
\nu^N + \frac{1}{2N} \Delta^N \nu^N 
= \left( h^N + \frac{1}{2N} \Delta^N h^N \right) 
[\Sigma^N]^{-1} (1^N,\nabla^N)^* \nu^N . 
\end{equation}
Since
\begin{equation}
\nu^N = 1^N \nu^N + N^{-1} \partial ^N_t \nu^N, \ 
h^N = 1^N h^N + N^{-1} \partial ^N_t h^N
\end{equation}
and
\begin{equation}
(1^N h^N)[\Sigma^N]^{-1} = (1,0,\ldots,0), 
\end{equation}
we have
\begin{equation}
\partial ^N_t \nu^N + \frac{1}{2} \Delta^N \nu^N = 
\left( \partial ^N_t h^N + \frac{1}{2} \Delta^N h^N \right) 
[\Sigma^N]^{-1} (1^N,\nabla^N)^* \nu^N . 
\end{equation}
This is exactly (\ref{disdif}). \qed
%%%

\subsection{Proof of Theorem \ref{limittheorem}}\label{PLMT}
The following proof is a routine-work in the context of
finite difference method.

First we will show that our scheme is {\em stable}.
Let $ u^N $ be the unique solution of the following difference equation.
\begin{equation}\label{disdif2}
\begin{split}
& u^N(T_N,x) = \Psi^N(x) ; \ x \in \bR^n, \\
& \partial ^N_t u^N + \frac{1}{2} \Delta^N u^N 
- \langle b^N,\nabla^N u^N \rangle - c^N (1^N u^N) = g^N ; \\
& \qquad t \in (0,T_N)_N , \ x \in \R^n.
\end{split}
\end{equation}
where $ g^N $ and $ \Psi^N $ are given functions on $ \R_+ \times \R^n $ and $ \R^n $ respectively. 
We claim that
\begin{equation}\label{stability}
\sup_{x \in \R^n} |u^N (t,x)| \leq \sup_{(s,y) \in [t,T] \times \R^n} \left\{
(T-t)|g^N (s,y)| + |\Psi^N (s,y)| \right\} 
\end{equation}
for every $ t \in [0,T]_N $. 
This inequality shows the stability of our scheme.
To prove (\ref{stability}), we first remark that the equation  in (\ref{disdif2}) can be rewritten
as
\begin{equation*}
1^N u^N = 1^N g^N /N + (1^N h^N) [H^N]^{-1} \tilde{u}^N, 
%\begin{split}
%&u^N (t,x) = g^N (t,x)/N + h^N (t,x)[H^N (t+ N^{-1},x)]^{-1} \\
%&\quad\times [u^N (t+N^{-1},x+N^{-1/2}e_1 ), ..., 
%u^N (t+ N^{-1},x +N^{-1/2}e_{n+1})]^*,
%\end{split}
\end{equation*}
which comes from Theorem \ref{Euro}.
By the positivity assumption on 
$ h^N(t,x) [H^N(t+N^{-1},x)]^{-1} $, 
we see that $ A^{-1} h^N(t,x) [H^N(t+N^{-1},x)]^{-1} $ 
defines a transition probability 
of a time-inhomogeneous Markov chain 
$ (Y^N_t,\bP^x_t)_{t \in [0,T]_N, x \in \bR^n} $: 
$ \bP^x_t(Y^N_{t}=x)=1 $ and
\begin{equation}
\begin{split}
& \bP^x_t(Y^N_{t+N^{-1}}=x+N^{-1/2}e_j) \\
& \mbox{$ = $ the $ j $-th component of} \\
& \quad A(t,x)^{-1} h^N(t,x) [H^N(t+N^{-1},x)]^{-1} . 
\end{split}
\end{equation}
Denoting the expectation with respect to 
$ \bP^x_t $ by $ \bE^x_t $, 
we have 
\begin{equation}\label{Markovian}
u^N(t,x) = \frac{1}{N} g^N(t,x) + 
\bE^x_t \left[ A(t,Y^N_t) u^N(t+N^{-1},Y^N_{t+N^{-1}}) \right]. 
\end{equation}
By iterating (\ref{Markovian}) and by the Markov property, we have 
\begin{equation}
\begin{split}
u^N(t,x) 
=& \frac{1}{N} 
\sum_{s \in [t,T]_N} \bE^x_t 
\left[ g^N(s,Y^N_{s}) \prod_{u \in [t,s)_N} 
A(u,Y^N_{u}) \right] \\
+& \bE^x_t \left[ 
\Psi^N(Y^N_{T_N}) 
\prod_{u \in [t,T_N)_N} A(u,Y^N_{u}) \right] . 
\end{split}
\end{equation}
(This is a discrete version of Feynman-Kac formula.)
By the assumption of $ 0< A < 1 $, we obtain (\ref{stability}).

Next, we will show that 
\begin{equation}\label{consistency}
\sup_N \sup_{(t,x) \in [0,T]_N \times \bR^n} 
N^{1/2} |g^N(t,x)| < \infty 
\end{equation}
where 
\begin{equation}\label{7.9}
g^N := \partial ^N_t \nu 
+ \frac{1}{2} \Delta^N \nu - \langle b^N,
\nabla^N \nu \rangle - c^N (1^N \nu)  
\end{equation}
for the solution $ \nu $ to (\ref{EuroPDE}). 
Since
\begin{equation}
g^N = g^N - \frac{\partial \nu}{\partial t} 
- \frac{1}{2} \Delta \nu + \langle b,\nabla \nu \rangle 
+ c \nu , 
\end{equation}
we have
\begin{equation}
\begin{split}
|g^N| \le& 
\left| \frac{\partial \nu}{\partial t} - \partial ^N_t \nu \right| 
+ \frac{1}{2} |\Delta \nu - \Delta^N \nu| \\
+& |b^N| |\nabla \nu - \nabla^N \nu| + |\nabla \nu| |b-b^N| \\
+& |c^N| |\nu - 1^N \nu| + |\nu| |c-c^N| . 
\end{split}
\end{equation}
By the Assumption \ref{standing} 
and the consistency (\ref{orders}), 
we obtain (\ref{consistency}).

Finally, by combining (\ref{stability}) and (\ref{consistency}),
and by the uniform continuity of $ \nu $, 
we have the desired result 
since $ \nu - \nu^N $ is the solution to (\ref{disdif2})
with $ \Psi^N(x) = \Phi(x)-\Phi^N(x) $ 
and $ g^N $ given by (\ref{7.9}). \qed

%%%%%%%%%%%%%%%%%%%%%%%%%%%%%%%%%%%%%%%%
%\begin{thebibliography}
\bibliographystyle{ims}
\bibliography{akahori+}

\def\polhk#1{\setbox0=\hbox{#1}{\ooalign{\hidewidth
  \lower1.5ex\hbox{`}\hidewidth\crcr\unhbox0}}} \def\cprime{$'$}
\begin{thebibliography}{12}
\expandafter\ifx\csname natexlab\endcsname\relax\def\natexlab#1{#1}\fi
\expandafter\ifx\csname url\endcsname\relax
  \def\url#1{\texttt{#1}}\fi
\expandafter\ifx\csname urlprefix\endcsname\relax\def\urlprefix{URL }\fi

\bibitem[{Akahori(2003)}]{Aka_Parisian}
\textsc{Akahori, J.} (2003).
\newblock Local time in parisian walkways.
\newblock Preprint, math.PR/0606183.

\bibitem[{Cox et~al.(1979)Cox, Ross and Rubinstein}]{CRR}
\textsc{Cox, J.}, \textsc{Ross, S.~A.} and \textsc{Rubinstein, M.} (1979).
\newblock Option pricing: a simple approach.
\newblock \textit{Journal of Financial Economics} \textbf{7} 229--263.

\bibitem[{Duffie(1996)}]{duffie:96}
\textsc{Duffie, D.} (1996).
\newblock \textit{Dynamic Asset Pricing Theory}.
\newblock 2nd ed. Princeton University Press, Princeton.

\bibitem[{Ethier and Kurtz(1986)}]{MR88a:60130}
\textsc{Ethier, S.~N.} and \textsc{Kurtz, T.~G.} (1986).
\newblock \textit{Markov Processes}.
\newblock Wiley Series in Probability and Mathematical Statistics: Probability
  and Mathematical Statistics, John Wiley \& Sons Inc., New York.

\bibitem[{Fujita(2002)}]{Fujita_Japanese_text}
\textsc{Fujita, T.} (2002).
\newblock \textit{Finance no tame no Kakuritu Kaiseki}.
\newblock Kodansha, Tokyo.
\newblock (in Japanese), [{\it Stochastic Caluculus for Finance.}].

\bibitem[{Fujita(2003)}]{Fujita_formula}
\textsc{Fujita, T.} (2003).
\newblock A random walk analogue of {L}\'evy's theorem.
\newblock Preprint, Hitotsubashi University.

\bibitem[{He(1990)}]{he:90}
\textsc{He, H.} (1990).
\newblock Convergence from discrete- to continuous-time contingent claims
  prices.
\newblock \textit{Review of Financial Studies} \textbf{3} 523--546.

\bibitem[{Kudzhma(1982)}]{MR684465}
\textsc{Kudzhma, R.} (1982).
\newblock It\^o's formula for a random walk.
\newblock \textit{Litovsk. Mat. Sb.} \textbf{22} 122--127.

\bibitem[{Protter(2004)}]{MR2020294}
\textsc{Protter, P.~E.} (2004).
\newblock \textit{Stochastic Integration and Differential Equations}, vol.~21
  of \textit{Applications of Mathematics (New York)}.
\newblock Springer-Verlag, Berlin.

\bibitem[{Richtmyer and Morton(1994)}]{MR95b:65003}
\textsc{Richtmyer, R.~D.} and \textsc{Morton, K.~W.} (1994).
\newblock \textit{Difference Methods for Initial-Value Problems}.
\newblock Robert E. Krieger Publishing Co. Inc., Malabar, FL.

\bibitem[{Serre(1978)}]{MR80f:20001}
\textsc{Serre, J.-P.} (1978).
\newblock \textit{Repr\'esentations lin\'eaires des groupes finis}.
\newblock Hermann, Paris.
\newblock [English translation {\em Linear Representations of Finite Groups.}
  Translated from the second French edition by Leonard L. Scott. Graduate Texts
  in Mathematics, Vol. 42. Springer-Verlag, New York-Heidelberg, 1977. ISBN
  0-387-90190-6].

\bibitem[{Szabados(1990)}]{MR92i:60105}
\textsc{Szabados, T.} (1990).
\newblock A discrete {I}t\^o's formula.
\newblock In \textit{Limit theorems in probability and statistics (P\'ecs,
  1989)}, vol.~57 of \textit{Colloq. Math. Soc. J\'anos Bolyai}. North-Holland,
  Amsterdam, 491--502.

\end{thebibliography}
%\end{thebibliography}
\end{document}